\documentclass{aims} 
\usepackage{amsmath,float}
\usepackage{paralist}
\usepackage[misc]{ifsym}
\usepackage{epsfig} 
\usepackage{epstopdf} 
\usepackage[colorlinks=true]{hyperref}
\hypersetup{urlcolor=blue, citecolor=red}
\allowdisplaybreaks
\def\theequation{\arabic{section}.\arabic{equation}}
\def\currentvolume{X}
 \def\currentissue{X}
  \def\currentyear{200X}
   \def\currentmonth{XX}
    \def\ppages{X-XX}
     \def\DOI{ }
\pdfpagesattr
{ /CropBox [55 40 546 730] }

\newtheorem{theorem}{Theorem}[section]

\newtheorem{lemma}[theorem]{Lemma}

\theoremstyle{definition}
\newtheorem{definition}[theorem]{Definition}

\newcommand{\eps}[1]{{#1}_{\varepsilon}}

\title[Weak solutions to  the time-fractional Fokker--Planck equation]
{Well-posedness and simulation of weak solutions to  the time-fractional Fokker--Planck equation with general forcing} 

\author[Marvin Fritz]{}


\thanks{The author is supported by the state of Upper Austria.}

\thanks{$^*$Corresponding author: Marvin Fritz}

\usepackage[T1]{fontenc}

\usepackage{xcolor,bm, amssymb, mathtools,subfigure}
\usepackage[capitalise,nameinlink]{cleveref}

\usepackage{tikz,bm}
\usepackage[sort,compress]{cite}

\newcommand{\pt}{\partial_t}
\renewcommand{\eps}{\varepsilon}
\newcommand{\Pta}{D_t^\alpha}
\newcommand{\pta}{\p_t^\alpha}

\newcommand{\ptb}{\p_t^{1-\alpha}}
\newcommand{\Ptb}{D_t^{1-\alpha}}

\newcommand{\ga}{g_{\alpha}}
\newcommand{\gb}{g_{1-\alpha}}

\newcommand{\ds}{\,\textup{ds}}

\newcommand{\dx}{\dd x}

\newcommand{\loc}{\textup{loc}}

\renewcommand{\O}{\mathcal{O}}

\renewcommand{\div}{\textup{div}}

\newcommand{\dd}{\mathop{}\!\mathrm{d}}
\renewcommand{\d}{\mathop{}\!\mathrm{d}}
\newcommand{\ddt}{\frac{\dd}{\dd\mathrm{t}}}
\newcommand{\dt}{\,\textup{d}t}
\newcommand{\p}{\partial}

\newcommand{\R}{\mathbb{R}} %

\renewcommand{\rho}{\varrho}

\newcommand{\con}{\hookrightarrow}

\renewcommand{\L}{\mathcal{L}} %
\newcommand{\W}{W} %

\newcommand{\longweak}{\relbar\joinrel\rightharpoonup}

\begin{document}
\maketitle

\centerline{\scshape
Marvin Fritz$^{{\href{mailto:marvin.fritz@ricam.oeaw.ac.at}{\textrm{\Letter}}}*1}$}

\medskip

{\footnotesize
 \centerline{$^1$Computational Methods for PDEs, Johann Radon Institute for Computational and}
 \centerline{Applied Mathematics, Linz, Austria}
} 

\bigskip



\begin{abstract}
In this paper, we investigate the well-posedness of weak solutions to the time-fractio\-nal Fokker--Planck equation. Its dynamics is governed by anomalous diffusion, and we consider the most general case of space-time dependent forces. Consequently, the fractional derivatives appear on the right-hand side of the equation, and they cannot be brought to the left-hand side, which would have been preferable from an analytical perspective.  
For showing the model's well-posedness, we derive an energy inequality by considering nonstandard and novel testing methods that involve a series of convolutions and integrations. We close the estimate by a Henry--Gronwall-type inequality. Lastly, we  propose a numerical algorithm based on a nonuniform L1 scheme and present some simulation results for various forces.
\end{abstract}


	\section{Introduction}
	

Mathematicians and engineers have given time-fractional differential equations significant consideration in recent years.
Such equations are nonlocal in time and possess an inherent history effect.
We refer the interested reader to the multi-volume work “Handbook of Fractional Calculus with Applications” \cite{baleanu2019handbook,petras2019handbook,tarasov2019handbook} and its references therein for more information and typical real-world applications such as physics, control theory, engineering, life and social sciences.

In this work, we are concerned with the time-fractional Fokker--Planck equation, which permits subdiffusive behavior and its derivation and application have been investigated earlier in literature. We distinguish the model by its exterior force, which can be time-dependent \cite{sokolov2006field,magdziarz2009stochastic}, space-dependent \cite{barkai2000continuous,barkai2001fractional,chechkin2003fractional,fu2021continuous,metzler1999anomalous,metzler1999deriving,metzler2002space,sandev2015diffusion,sokolov2001dynamics},  or space-time dependent \cite{angstmann2015generalized,heinsalu2007use,magdziarz2008equivalence,weron2008modeling}. We focus on the latter, most general, case and mention the publications \cite{huang2020new,le2016numerical,le2018semidiscrete,le2021alpha,mclean2021uniform,mustapha2022second,pinto2017numerical,yan2019finite} that explored numerical methods for the time-fractional Fokker--Planck equation with space-time dependent forces. 

We emphasize that several articles have investigated a “time-fractional Fokker--Planck”-type equation, where the time-fractional derivative in the sense of Caputo appears on the left-hand side of the PDE. This is correct in the case of a time-independent force. However, for time dependent forces this model is not correct and according to  \cite{heinsalu2007use}, it is “physically defeasible” and its solution “does not correspond to a physical stochastic process”. In this work, we provide some mathematical and numerical insights on this reformulation and the differences of both models.

We present an analytical treatment of weak solutions to the time-fractional Fokker--Planck equation with space-time dependent forces. Specifically, we follow the Galerkin ansatz by spatially discretizing the system and deriving appropriate energy constraints, allowing us to reach the limit in the discretized system.
We mention that weak solutions to other nonlinear time-fractional PDEs have been previously studied using the Galerkin method in the published works \cite{fritz2021subdiffusive,fritz2023equivalence,fritz2022time}.
In addition, preliminary steps have been taken in the optimal control \cite{camilli2020approximation} and analysis \cite{mclean2020regularity,le2019existence,le2021alpha,mclean2021uniform,peng2022existence,mclean2019well} of the time-fractional Fokker–Planck system.
Nonetheless, mild, strong, and classical solutions have been investigated. The difficulty lies in the low regularity of  weak solutions and the appearance of time-dependent forces, which do not allow us to transform the system to a more accessible system regarding analysis. A coupled system of a time-fractional Fokker--Planck equation with the Navier--Stokes equations was investigated in the work \cite{fritz2023analysis} but because of the complex coupling between the equations only the case of $\alpha \in (\frac12,1)$ was considered.

In Section \ref{Sec:Derivation}, we discuss the mathematical model with its initial and boundary data.
In Section \ref{Sec:Prelim}, we present various function spaces and recall important conclusions from the theory of fractional derivatives, including chain inequalities, embedding theorems, and Gronwall-type inequalities.
In Section \ref{Sec:Analysis}, we finally present and verify the theorem declaring the well-posedness of weak solutions.
Here, the system is discretized, and appropriate energy bounds are derived to pass the limit in the discretized system. 
In Section \ref{Sec:Numerics}, we propose a numerical discretization of the time-fractional equation based on the nonuniform L1 scheme in time and finite elements in space. We show simulation results and focus on the influence of the fractional derivative. Moreover, we compare the model that is studied here to the physical defeasible model as mentioned above.

	\section{Modeling of the time-fractional Fokker--Planck equation} \label{Sec:Derivation}


 Let $\Omega \subset \R^d$, $d\in \mathbb{N}$, be a Lipschitz domain and $T<\infty$ a fixed final time. Shortly, we denote the time-space domain by $\Omega_T=\Omega \times (0,T)$. Let $\psi:\Omega_T \to \R$ denote a probability density function that represents the probability at a time $t$ of finding the center of mass of a particle in the volume element $x+\d x$.

The time-fractional Fokker--Planck model with space-time dependent force can be derived by utilizing the Langevin equations, see \cite{magdziarz2008equivalence,magdziarz2009stochastic}, and the model reads
\begin{equation} \label{Eq:DerivFP}
\pt \psi(x,t)- D \Delta \Ptb \psi(x,t)+\div\!\left(F(x,t) \Ptb \psi(x,t) \right)  = 0.\end{equation}
Here, $F:\Omega_T  \to \R^d$ denotes the space-time dependent external force and $D$ the diffusion coefficient. In contrast to the typical model of integer-order, the fractional derivative in the sense of Riemann--Liouville is introduced, which is defined by
$$\Ptb u(t) =\frac{1}{\Gamma(\alpha)} \ddt \int_0^t \frac{u(s)}{(t-s)^{1-\alpha}} \ds, $$
where $\Gamma$ denotes Euler's Gamma function. We introduce the singular kernel function $g_\alpha(t)=t^{\alpha-1}/\Gamma(\alpha)$ and therefore, we can rewrite the fractional derivative with the convolution operator as $$\Ptb u =\pt (\ga*u).$$ In the limit case of $\alpha=1$, the model is reduced to the standard Fokker--Planck equation.
This time-fractional model has been studied in the previous works \cite{huang2020new,le2016numerical,le2018semidiscrete,le2021alpha,mclean2021uniform,mustapha2022second,pinto2017numerical,yan2019finite} with regard to numerical methods and in \cite{le2019existence} for the existence of mild and classical solutions. 

We note that the fractional derivative in the sense of Riemann--Liouville appears naturally in the equation's derivation, see \cite{magdziarz2008equivalence}. However, the fractional derivative in the sense of Caputo would be preferable considering our variational approach to time-fractional partial differential equations and the involved analytical machinery. The Caputo derivative of order $\alpha$ is denoted by $\pta$ and it reads 
\begin{equation} \label{Eq:Caputo} \pta u= \Pta (u-u_0).
\end{equation} 
Here, $u_0$ is the initial of the underlying system, which shall fulfill $$\big(\gb*(u-u_0)\big)(0)=0$$ in the case that $u$ is not continuous. 

If the force is time-independent, we could simply convolve the time-fractional Fokker--Planck equation \eqref{Eq:DerivFP} with the singular kernel function $\gb$ and exploit the properties $\gb*\Ptb u=u$ and $\pta u=\gb*\pt u$, see below in Section \ref{Sec:Prelim}, to obtain the time-fractional equation
\begin{equation} \label{Eq:DerivFP3}
\pta \psi(x,t)-D \Delta \psi(x,t) +\div \big(F(x) \psi(x,t) \big)  =0,\end{equation}
which would be more accessible for analytical and numerical methods. However, we cannot simply exclude the relevant cases of time-dependent forces. In such cases, one would require a product rule for fractional derivatives to write $F\Ptb \psi$ as $\Ptb(F\psi)-\Ptb F \psi$.  However, this is not correct for fractional derivatives, as it can be already seen from the example $\psi=F=1$. Then it holds 
$$F\Ptb \psi = \ga \neq 0 = \ga-\ga =\Ptb(F\psi)-\Ptb F \psi.$$
There is a fractional version of the Leibniz rule that requires two smooth functions $f,g$ and reads \cite[Theorem 2.18]{diethelm2010analysis}
$$\Pta(fg)=f \Pta g + \sum_{k=1}^\infty \binom{\alpha}{k} \pt^k f \cdot (g_{1-k+\alpha}*g).$$
We can already see the issue of this formula. It requires smooth functions, and it turns out that there is an infinite sum on the right-hand side. Let us assume that $F$ and $\psi$ are smooth. Then we want to bring the fractional derivative in front of $F\psi$ by the formula
$$F \Ptb \psi = \Ptb(F\psi) - \sum_{k=1}^\infty \binom{1-\alpha}{k} \pt^k F \cdot (g_{2-k-\alpha}*\psi).$$
Afterward, we convolve the system with $\gb$ and obtain the system
\begin{equation} \label{Eq:ModelWrong} \begin{aligned}
&\pta \psi(x,t)-D \Delta \psi(x,t) +\div \big(F(t,x) \psi(x,t) \big)  \\&=\sum_{k=1}^\infty \binom{1-\alpha}{k} \gb*\big(\pt^k F \cdot (g_{2-k-\alpha}*\psi)\big),
\end{aligned}\end{equation}
There have been several published articles that studied this model but neglecting the complete right-hand side. This is also the reason it is claimed in \cite{heinsalu2007use} that such a model (with neglecting the right-hand side) is “physically defeasible” and its solution “does not correspond to a physical stochastic process”.  In the case that $F$ is affine linear in $t$, i.e. $F(t,x)=a(x)+b(x)t$,  it yields 
\begin{equation*} \begin{aligned}
&\pta \psi(x,t)-D \Delta \psi(x,t) +\div \big(F(t,x) \psi(x,t) \big)  \\&=(1-\alpha) \cdot b(x) \cdot (g_{2-2\alpha}* \psi)(t)
\end{aligned}\end{equation*}
 We would rather not consider infinitely many terms on the right-hand side of the PDE for a general $F$ and therefore, we instead 
exploit the definition \eqref{Eq:Caputo} of the Caputo derivative to obtain
$$\Ptb u(t) = \ptb u(t) + \Ptb u_0=\ptb u(t)+u_0g_{\alpha}(t),$$
and rewrite the time-fractional Fokker--Planck equation \eqref{Eq:DerivFP} as follows:
\begin{equation} \label{Eq:DerivFP2}
\begin{aligned}
&\pt \psi(x,t)-D \Delta \ptb \psi(x,t) +\div\big(F(x,t) \ptb \psi(x,t) \big) \\ &=\ga D\Delta  \psi_0 - \ga \div(F\psi_0).
\end{aligned} \end{equation}
We consider an initial condition $\psi_0 \in H_0^1(\Omega)$ and therefore, it holds that the right-hand side has the regularity $L^p(0,T;H^{-1}(\Omega))$ with $p<1/(1-\alpha).$
We equip this equation with the homogeneous Dirichlet boundary condition $\psi=0$ on $\p\Omega$. However, our analytical results also hold for no-flux boundary conditions (i.e. homogeneous Neumann). Moreover, the system is equipped
with the initial condition $\psi(0) = \psi^0 \geq 0$ in  $\Omega$. 
Physically, $\psi^0$ is a given probability density function, i.e.,  it is nonnegative function and satisfies $\int_\Omega \psi^0(x) \dx =1$ (however, we do not need to assume such properties in our well-posedness theorem below).   
Integrating the time-fractional Fokker--Planck equation  in $\Omega$  and employing integration by parts, we find
$\ddt \int_{\Omega} \psi(x,t) \dx =0.$ This implies then $\int_{\Omega} \psi(x,t) \dx = 1$  for almost all $t$.

\section{Mathematical preliminaries} \label{Sec:Prelim}
In this part, we present some important concepts and conclusions addressing fractional derivatives.
For instance, we provide a fractional version of the Aubin--Lions lemma and a suitable Gronwall lemma. These are important results used in Galerkin-based proofs for showing the existence of weak solutions to partial differential equations.


Let $T<\infty$ be a fixed final time. We have already defined the singular kernel function in the previous section by $g_\alpha(t)=t^{\alpha-1}/\Gamma(\alpha)$, $t \in (0,T)$, $\alpha >0$. We can extend the definition to the limit case of $\alpha=0$ by $g_0=\delta$. We observe that it holds  $g_\alpha \in L^p(0,T)$ for any $\alpha>1-\frac{1}{p}$, i.e., 
\begin{equation}\label{Eq:GaLp}
\ga \in L^{\frac{1}{1-\alpha}-\eps}(0,T) \quad \forall \eps \in \big(0,\tfrac{\alpha}{1-\alpha}\big].
\end{equation} 
Alternatively, using the concept of locally integrable functions, it naturally holds $\ga \in L_\loc^{1/(1-\alpha)}(0,T)$.
E.g., it holds $\ga \in L^2(0,T)$ for any $\alpha>\frac12$ and $\ga \in L_\loc^2(0,T)$ for any $\alpha \geq \frac12$. Moreover, the kernel function satisfies the following semigroup property, see \cite[Theorem 2.4]{diethelm2010analysis},
\begin{equation} \label{Eq:Semigroup}
	\ga * g_\beta = g_{\alpha+\beta} \qquad \forall \alpha,\beta \in (0,1).
\end{equation} 
We note that one can bound the $L^p(0,t)$-norm of a function $u:(0,T) \to \R$ by its convolution with $\ga$ as follows:
\begin{equation}\begin{aligned} \|u\|_{L^p_t}^p := \int_0^t |u(s)|^p \ds  &\leq t^{1-\alpha} \int_0^t (t-s)^{\alpha-1} |u(s)|^p \ds \\ &\leq T^{1-\alpha} \Gamma(\alpha) \big(\ga * |u|^p\big)(t).	
\end{aligned} 
\label{Eq:KernelNorm}
\end{equation}
In other words, the space
\begin{equation} \label{Eq:LpAlpha} L^p_\alpha(0,T)=\Big\{u:(0,T) \to \R: \|u\|_{L^p_\alpha}^p:=\sup_{t \in (0,T)}(\ga*|u|^p)(t) < \infty\Big\},
\end{equation}
is indeed continuously embedded in the space $L^p(0,T)$.
We can relate the estimate  \eqref{Eq:KernelNorm} to $\ga$ by noting that
\begin{equation*} \begin{aligned} (\ga* |u|^p)(t) \geq \frac{t^{\alpha-1}}{\Gamma(\alpha)} \|u\|^p_{L^p_t} = \ga(t) \|u\|^p_{L^p_t} \geq \ga(T) \|u\|^p_{L^p_t}.
\end{aligned}
\end{equation*}
In particular, this yields for any $s \leq t$
\begin{equation} \label{Eq:EstimateST} (\ga* |u|^p)(t) \geq (\ga* |u|^p)(s) \geq  \ga(s) \|u\|^p_{L^p_s}.
\end{equation}
Therefore, we can integrate this inequality on the time interval $(0,t)$ to obtain
$$t \cdot (\ga* |u|^p)(t)  \geq  \int_0^t  \ga(s) \|u\|^p_{L^p_s} \ds,$$
which implies the following useful bound
\begin{equation}\label{Eq:IneqGaG1}\begin{aligned}(\ga* |u|^p)(t)  &\geq \frac{1}{t} \int_0^t  \ga(s) \|u\|^p_{L^p_s} \ds \\ &\geq  \frac{1}{T} \int_0^t  \ga(s) \|u\|^p_{L^p_s} \ds .\end{aligned}\end{equation}
Similarly, if we take the convoluton instead of the integration of the inequality \eqref{Eq:EstimateST}, we obtain
\begin{equation}\label{Eq:IneqGaG1Conv}\begin{aligned} g_{\alpha+1}(T) (\ga* |u|^p)(t) &\geq g_{\alpha+1}(t) (\ga* |u|^p)(t) \\  &\geq  \big(\ga * (\ga\cdot \|u\|^p_{L^p_t})\big)(t).
\end{aligned}
\end{equation}


In the previous section, have also introduced the fractional derivatives in the sense of Riemann--Liouville $\Pta u = \pt(\gb*u)$ and Caputo $\pta u=\Pta(u-u_0)$. It is well-known that the Caputo derivative can also be written as $\pta u=\gb*\pt u$ if $u$ is absolutely continuous, see \cite[Lemma 3.5]{diethelm2010analysis}. We note that it does not hold $\pta \pt^\beta u = \pt^{\alpha+\beta} u$ in general for the Caputo derivative. However, it holds, see \cite[Theorem 3.14]{diethelm2010analysis}, 
\begin{equation} \label{Eq:PtbPtaPt}
\pta\ptb u = \pt u.
\end{equation}
We define the fractional Sobolev--Bochner space   for $\alpha \in (0,1)$ on $(0,T)$ with values in a given Hilbert space $H$ by $$\W^{\alpha,p}(0,T;H)=\big\{u \in L^p(0,T;H) : \pta u \in L^{p}(0,T;H)\big\}.$$
Next, we state the inverse convolution property. Its name origins from the fact that the convolution with the kernel $\ga$ acts as an inverse operation on the $\alpha$-th fractional derivative up to the initial condition. In fact,  it holds
\begin{align} \label{Eq:InverseConvolution}
		(\ga* \pta u)(t)    &=u(t)-u_0  \qquad \forall u \in \W^{\alpha,p}(0,T;H). \end{align} 
  This can be seen from the computation
  $$(\ga*\pta u)(t)=(\ga*\gb *\pt u)(t)=(1*\pt u)(t)=\int_0^t \pt u(s) \ds = u(t)-u_0,$$
  where we used \eqref{Eq:Semigroup} to conclude $\ga*\gb=g_1=1$.
Furthermore, we mention the following consequences of the interaction between fractional derivatives and kernel functions:
\begin{equation} \label{Eq:DerivativeofKernel}  \begin{aligned}
	\pta (\ga * u ) &=\Pta (\ga * u)=\pt ( \gb * \ga * u) = \pt (1*u) = u,
 \end{aligned}
\end{equation}
which holds for any $u \in L^1(0,T;H)$. 




As in the integer-order setting, there are continuous and compact embedding results for fractional Sobolev spaces; see \cite[Theorem 3.2]{wittbold2021bounded}. For a given Gelfand triple $V \con\con H \con V'$, the classical Aubin--Lions lemma \cite{simon1986compact} reads
\begin{equation}\label{Eq:aubin}  \begin{aligned}
\W^{1,1}(0,T;V') \cap L^\infty(0,T;V) &\con\con C([0,T];H),\\ %
\W^{1,1}(0,T;V') \cap L^p(0,T;V) &\con\con L^p(0,T;H), \quad p \in [1,\infty),
\end{aligned}\end{equation} 
and the fractional counterparts is as follows:
\begin{equation}\label{Eq:aubinfractional}  \begin{aligned}
\W^{\alpha,p}(0,T;V') \cap L^{p'}(0,T;V) &\con C([0,T];H), &&p \in [1,\infty),\\%
\W^{\alpha,p}(0,T;V') \cap L^p(0,T;V) &\con\con L^p(0,T;H), &&p \in [1,\infty).
\end{aligned}\end{equation} 
We observe that there is a give-and-take involved: The fractional derivative is of order $\alpha<1$ i.e. it is less than the full derivative in the classical Aubin--Lions lemma. However, we require that the derivative is in the better space $L^p(0,T;V')$ instead of only $L^1(0,T;V')$ to achieve the same target space $L^p(0,T;H)$ in the compactness result.

The classical chain rule does not hold for fractional derivatives, but one can use the following inequality, see \cite[Theorem 2.1]{vergara2008lyapunov}, as a remedy:
\begin{equation} \label{Eq:ChainOriginal}  \frac12 \pta \|u\|^2_H  \leq (u,\pta u)_H \quad \forall u \in W^{\alpha,p}(0,T;H),
\end{equation}
for almost all $t \in (0,T)$, which is also known as Alikhanov's inequality, see the original work \cite{alikhanov2010priori}. Moreover, we conclude from \eqref{Eq:PtbPtaPt} and Alikhanov's inequality the following:
\begin{equation} \label{Eq:ChainExtended} (\pt u, \pta u)_H = (\ptb \pta u,\pta u)_H \geq \frac12 \ptb \|\pta u\|_H^2,
\end{equation}
which gives after integrating it over the time interval $(0,t)$
$$\int_0^t (\pt u, \pta u)_H \ds \geq \frac12 (\ga * \|\pta u\|_H^2)(t) \geq \frac{1}{2\Gamma(\alpha)T^{1-\alpha}} \|\pta u\|_{L^2_tH}^2,
$$
where we applied \eqref{Eq:KernelNorm} in the last step.


Next, we require a Gronwall-type inequality that allows convolutions on the right-hand side of the inequality. Moreover, we want to have an additional function on the right-hand side that is only  locally integrable. Such inequalities are known as Henry--Gronwall inequalities.

\begin{lemma}[{Henry--Gronwall, cf. \cite[Lemma 7.1.1]{henry1981geometric}}]
    Let $b\geq 0$, $\beta>0$, $a \in L^1_\loc(0,T;\R_{\geq 0})$. If $u \in L^1_\loc(0,T;\R_{\geq 0})$ satisfies
    $$u(t) \leq a(t)+b (g_\beta * u)(t), \quad \text{ for a.e. } t \in (0,T),$$
    then it yields
    $$u(t) \leq C(\alpha,b,T) \cdot \big((g_0+E)*a\big)(t) , \quad \text{ for a.e. } t \in (0,T),$$
    where $E$ is related to the Mittag--Leffler function.
\end{lemma}

We prove the following extension of the Henry--Gronwall inequality that allows an additional term on the left-hand side.

\begin{lemma}\label{Lem:GronFrac}
    Let $b\geq 0$, $b>0$, $a \in L^1_\loc(0,T;\R_{\geq 0})$. If the  functions $u,v \in L_\loc^1(0,T;\R_{\geq 0})$  satisfy the inequality
    $$u(t) + (\ga * v)(t) \leq a(t)+b (\ga* u)(t) \qquad \text{for a.a. } t \in (0,T], $$
    then it yields
    $$u(t) + \int_0^t v(s) \ds \leq C(\alpha,b,T) \cdot \big( (g_0+E)*a\big)(t) \qquad \text{for a.a. } t \in (0,T]. $$
\end{lemma}
\begin{proof}
We define the function $w=u+\ga*v$.
    Since $\ga*v$ is again nonnegative, we obtain
    $$w(t) \leq a(t)+b (\ga* u)(t) \leq a(t) + (\ga*w)(t),$$
    and by the Henry--Gronwall it yields
    $$w(t) \leq C(\alpha,b,T) \cdot \big( (g_0+E)*a\big)(t).$$
    Moreover, we can use \eqref{Eq:KernelNorm} to estimate $\ga*v$ by the integral of $v$, and we obtain the lemma's desired bound.
\end{proof}

\section{Well-posedness of weak solutions} \label{Sec:Analysis}
In this section, we state and prove the well-posedness of weak solutions to the time-fractional Fokker--Planck equation \eqref{Eq:FP}. As we already mentioned, we equip the equation with a homogeneous Dirichlet boundary condition. As noted before, our analysis holds for no-flux boundary conditions as well. 
We analyze the PDE in the Hilbert triple 
$$H_0^1(\Omega) \hookrightarrow \hookrightarrow L^2(\Omega) \hookrightarrow H^{-1}(\Omega).$$
We equip $H_0^1(\Omega)$ with the norm $\|\cdot\|_{H_0^1}=\|\nabla \cdot \|_{L^2(\Omega)}$ with is equivalent to the natural norm on $H^1(\Omega)$ due to Poincar\'e's inequality \cite[6.7]{alt2016linear}.
We use the Galerkin method and discretize the partial differential equations in space. Further, we derive suitable energy estimates, and we emphasize the places where the time-fractional derivative comes into play. We shall then pass to the limit to deduce the existence of a weak solution. The uniqueness is obtained as usual.

First off, however, we introduce the concept of a  weak solution to the time-fractional Fokker--Planck equation in the following definition.\medskip

\begin{definition} \label{Def:Weak}
We call a function $\psi:\Omega_T \to \R$ a weak solution to the time-fractional Fokker--Planck equation \eqref{Eq:FP} if it is of the regularity
$$\psi \in W^{1,1}(0,T;H^{-1}(\Omega)) \cap H^{1-\alpha}(0,T;H_0^1(\Omega)),$$
fulfills the initial data  $\psi(0)=\psi_0$ in $H^{-1}(\Omega)$, and the following variational form: 
\begin{equation}\label{Eq:FP} \begin{aligned}
  &\langle \pt \psi,\zeta \rangle_{H_0^1}
  +  D(\ptb \nabla \psi,\nabla \zeta)_{L^2} - (F \ptb \psi,\nabla \zeta)_{L^2} \\ &=\langle f,\zeta\rangle_{H_0^1} - \ga (D \nabla \psi_0,\nabla \zeta)_{L^2} + \ga \cdot (F \psi_0,\nabla \zeta)_{L^2} \qquad \forall \zeta \in H_0^1(\Omega).
\end{aligned} \end{equation}
\end{definition} \medskip

As we see, we expect a solution that is continuous in time with values in the Hilbert space $H^{-1}(\Omega)$. Therefore, it is well-defined for the initial to fulfill $\psi(0)=\psi_0$ in $H^{-1}(\Omega)$. Moreover, it holds
$$H^{1-\alpha}(0,T;H_0^1(\Omega)) \hookrightarrow C([0,T];H_0^1(\Omega)),$$
if $1-\alpha>1/2$ i.e. $\alpha<\frac12$. In this case, the initial is even satisfied in $H_0^1(\Omega)$. 

Next, we state the main result of this work on the well-posedness of weak solutions to the time-fractional Fokker--Planck equation \eqref{Eq:FP}.

\begin{theorem}[Well-posedness of weak solutions] \label{Thm:WellPosedness} Let us assume: 
\begin{itemize}
	\item $\Omega \subseteq \R^d$, $d \in \mathbb{N}$, bounded Lipschitz domain,  $T<\infty$ fixed final time,
 \item $\alpha \in (0,1)$,
       \item $\psi_0 \in H_0^1(\Omega)$,
 \item $f \in L^2_\alpha(0,T;H^{-1}(\Omega))$,
\item $F \in L^\infty(\Omega_T;\R^d)$ with $\|F\|_{L^\infty(\Omega_T)} \leq F_\infty<\infty$.
  \end{itemize}
	Then there exists a unique weak solution $\psi$ to the time-fractional Fokker--Planck equation \eqref{Eq:FP} in the sense of Definition \ref{Def:Weak}. Further, it has the additional regularity
 $$\psi \in W^{1,r'}(0,T;H^{-1}(\Omega)) \cap W^{1-\alpha,p}(0,T;L^2(\Omega)) \cap H^{1-\alpha}(0,T;H_0^1(\Omega)),$$
 with $r'$ being the H\"older conjugate of $r=\max\{q',2\}$, $q'$ being the H\"older conjugate of $q=\frac{1}{1-\alpha}-\eps$ for $\eps\in (0,\frac{\alpha}{1-\alpha}]$, and $$\begin{cases}
	p=\infty, &\alpha>\frac12, \\
	p<\infty, &\alpha=\frac12,\\
	p<\frac{2}{1-2\alpha}, &\alpha<\frac12.
\end{cases} $$
\end{theorem}

We comment on the assumptions in this well-posedness result. We see it as an advantage that we can show the equation's well-posedness for all fractional values between $0$ and $1$, and any dimension $d \geq 1$. Further, we only require $F \in L^\infty(\Omega_T;\R^d)$ as opposed to \cite{le2019existence} that required $F \in W^{2,\infty}(\Omega_T)$ for showing results on mild solutions. Moreover, the work \cite{mclean2019well} studied a Volterra integral form of a class of time-fractional advection-diffusion-reaction equations, including the time-fractional Fokker--Planck equations. However, they required $F \in C^2([0,T];W^{1,\infty}(\Omega)^d)$ to show the well-posedness of the Volterra integral equation.

As the solution lies in the space $W^{1-\alpha,\infty}(0,T;L^2(\Omega))$ for $\alpha>\frac12$, we obtain $\psi \in C([0,T];L^2(\Omega))$ for $\alpha>\frac12$. It remains to consider $\alpha=\frac12$. In this case, we have $q=2-\eps$ and $r=\max\{1-1/(2-\eps),2\}=2$. Therefore, it holds $\psi \in H^1(0,T;H^{-1}(\Omega)) \cap L^2(0,T;H_0^1(\Omega))$, i.e., by an interpolation result $\psi \in C([0,T];L^2(\Omega))$. We summarize the continuity results as follows:
$$\psi \in \begin{cases}
    C([0,T];L^2(\Omega)), &\alpha \in (0,1), \\
    C([0,T];H_0^1(\Omega)), &\alpha \in (0,\frac12).
\end{cases}$$
We conclude that the initial is indeed at least satisfied in $L^2(\Omega)$. 

\bigskip




\noindent\textbf{Proof} In order to prove this theorem, we  employ the Galerkin method to discretize the variational form in space. This reduces the time-fractional PDE to a system of fractional ODEs, which admits a discretized solution $\psi_k$. We then derive $k$-uniform energy estimates, which imply the existence of weakly/weakly-$*$ convergent subsequence $\psi_{k_j}$. Finally, we pass to the limit $j \to \infty$ and apply compactness methods to return to the variational form of the continuous system. Recently, the Galerkin method has been applied to various time-fractional PDEs, see, e.g., \cite{fritz2021subdiffusive,fritz2022time,fritz2023equivalence,vergara2015optimal}. \medskip

\noindent \textbf{(1) Galerkin discretization.}
We introduce the discrete spaces
\begin{align*}
	H_k & =\text{span}\{ h_1,\dots,h_k\}, \\
\end{align*}
where $h_j: \Omega \to \R$, 
$j \in \{1,\dots,k\}$, are the eigenfunctions to the eigenvalues $\lambda_{j} \in \R$ of the following problems
$$\begin{aligned}
	(\nabla h_j,\nabla v)_{L^2} &= \lambda_{j} (h_j,v)_{L^2} &&\forall  v \in H_0^1(\Omega).  
\end{aligned}$$

Since the inverse Dirichlet--Laplace operator is compact, self-adjoint, injective, positive operators on $L^2(\Omega)$, we conclude by the spectral theorem, see e.g.,\cite[12.12 and 12.13]{alt2016linear}, that
\begin{alignat*}{3}
	& \{h_j\}_{j \in \mathbb{N}} &  & \text{ is an orthonormal basis in } L^2(\Omega) &  & \text{ and orthogonal in } H_0^1(\Omega), 
\end{alignat*}
Therefore, ${\cup_{k\in\mathbb{N}}} H_k$ 
is dense in $H_0^1(\Omega)$. We consider the Galerkin approximations
\begin{equation}\begin{gathered}
		\psi_k (t) = \sum_{j=1}^k \psi^j_k(t) y_j,
	\end{gathered}
	\label{Eq:GalerkinAnsatzFunctions}
\end{equation}
where 
$\psi^j_k: (0,T) \to \R$
are coefficient functions for all $j \in \{1,\dots,k\}$.
We denote the orthogonal projections onto the finite-dimensional space by $\Pi_{H_k}: L^2(\Omega) \to H_k$. 
Given the initial data $\psi_0$ from the continuous system, we choose $\psi_{0k} \in H_k$ such that  $\psi_{0k}=\Pi_{H_k} \psi_0$, i.e., there are coefficient $\{\psi_{0k}^{j}\}_{j=1}^k$ such that $\psi_{0k}=\sum_{j=1}^k \psi_{0k}^{j} y_j$. Moreover, due to well-known properties of the projection operator, see \cite[9.7]{alt2016linear}, it holds as $k \to \infty$
\begin{equation} \label{Eq:Projection}
   \|\psi_{0k}\|_X \leq \|\psi_0\|_X ~\text{ and }~ \psi_{0k} \to \psi_0 ~\text{ in }~ X \in \{H^{-1}(\Omega), L^2(\Omega), H_0^1(\Omega) \}.
\end{equation}

The Galerkin equations read as follows: We want to find  $\psi_k \in H_k$ such that $\psi_k(0)=\psi_{0k}$ and
\begin{equation}\label{Eq:FP_dis}\begin{aligned}
	&(\pt\psi_k,\zeta )_{L^2}
  +  D(\ptb \nabla \psi_k,\nabla \zeta)_{L^2} - (F \ptb \psi_k,\nabla \zeta)_{L^2} \\ &=\langle f,\zeta\rangle_{H^1_0} - \ga(D\nabla \psi_{0k},\nabla \zeta)_{L^2} + \ga \cdot (F\psi_{0k},\nabla \zeta)_{L^2}. 
 \end{aligned}\end{equation} for all $\zeta \in H_k$. We want to apply an existence result on ODEs with Riemann--Liouville derivatives and therefore, we rewrite the Galerkin system as follows:
 \begin{equation*}\begin{aligned}
	&(\pt\psi_k,\zeta )_{L^2}
  +  D(\Ptb \nabla \psi_k,\nabla \zeta)_{L^2} - (F \Ptb \psi_k,\nabla \zeta)_{L^2} =\langle f,\zeta\rangle_{H^1_0}, 
 \end{aligned}\end{equation*}
 for any $\zeta \in H_k$. We rewrite $\psi_k$ as the sum of the basis functions $\{\psi_k^j\}_{j=1}^k$ as introduced in \eqref{Eq:GalerkinAnsatzFunctions}, from which we obtain that the coefficients
are governed by the system
 \begin{equation}\label{Eq:FP_dis3}\begin{aligned}
	&\pt\psi_k^i 
  +  \lambda_i D\Ptb\psi_k^i - \sum_{j=1}^k \Ptb \psi_k^j (F y_j,\nabla y_i)_{L^2}=\langle f,y_i\rangle_{H^1},  
 \end{aligned}\end{equation}
 for any $i \in \{ 1, \dots, k\}$.
 Equivalently, we define the function $\phi_k^i=\Ptb\psi_k^i$ for any $i$, which is governed by the equation
 \begin{equation}\label{Eq:FP_dis2}\begin{aligned}
	&\Pta\phi_k^i 
  +  \lambda_i D \phi_k^i - \sum_{j=1}^k \phi_k^j (F y_j,\nabla y_i)_{L^2}  =\langle f,y_i\rangle_{H^1}, 
 \end{aligned}\end{equation} for any $i \in \{ 1, \dots, k\}$. We notice that it holds $\gb*\phi_k^i=\psi_k^i$ due to the inverse convolution property \eqref{Eq:InverseConvolution} and we observe that \eqref{Eq:FP_dis2} is naturally equipped with the initial $(\gb*\phi_k^i)=\psi_{0k}^i$.

 We denote the vector of components $\big(\phi_k^j(t)\big)_{1\leq j\leq k}$
by $\Phi(t)$. Then the approximate problem can be written as a system of ordinary
differential equations for $\Phi(t)$ of the form
$\Pta \Phi  = h(t,\Phi)$,
where $h$ is continuous and locally Lipschitz continuous with respect to $\Phi$. Therefore, the  fractional variant of the Cauchy--Lipschitz theorem, see \cite[Theorem 5.1]{diethelm2010analysis},  yields the existence of a unique continuous solution, defined on a short-time interval $[0, T_k]$ with $0<T_k \leq T$. From here, we conclude $\phi_k + \ga \psi_{0k} = \ptb \psi_k \in C((0,T_k];H_k)$ and $\ga*\psi_k \in C^1((0,T_k];H_k)$.
\medskip

\noindent \textbf{(2) Energy estimates: Part 1.} Next, we derive $k$-uniform estimates that will allow us to extract weakly converging subsequences. 
We test the Galerkin equation \eqref{Eq:FP_dis} by $\ga*\pt \psi_k=\ptb \psi_k \in H_k$ giving
\begin{equation} \label{Eq:Tested1}\begin{aligned} &(\pt \psi_k,\ptb \psi_k)_{L^2} + D\|\nabla \ptb \psi_k\|^2_{L^2} -  (F\cdot \nabla\ptb \psi_k,\ptb \psi_k)_{L^2} \\  &=  \langle f,\ptb \psi_k\rangle_{H_0^1} - \ga( D\nabla \psi_{0k} - F\psi_{0k},\nabla\ptb \psi_k)_{L^2}
\end{aligned}\end{equation}
For the first term on the left-hand side, we use $\pt \psi_k=\pta \ptb \psi_k$, see \eqref{Eq:PtbPtaPt}, to conclude with Alikhanov's inequality, see \eqref{Eq:ChainExtended},
\begin{equation} \label{Eq:Tested1a}
\begin{aligned} (\pt \psi_k,\ptb \psi_k)_{L^2} &=(\pta \ptb \psi_k,\ptb \psi_k)_{L^2} \\ &\geq \frac12 \pta \|\ptb \psi_k\|_{L^2}^2.
\end{aligned}
\end{equation}
We bring the term involving the force $F$ to the right-hand side of \eqref{Eq:Tested1} and apply the H\"older inequality to conclude
$$(F\nabla\ptb \psi_k,\ptb \psi_k)_{L^2} \leq F_\infty \|\nabla\ptb \psi_k\|_{L^2} \|\ptb \psi_k\|_{L^2},$$
where $F_\infty<\infty$ is the constant as introduced in the theorem's assumptions.
Further, we apply the Young inequality to give the norm of $\nabla \ptb \psi_k$ a prefactor that is smaller than $D$, i.e., we obtain
\begin{equation} \label{Eq:Tested1b}F_\infty \|\nabla\ptb \psi_k\|_{L^2} \|\ptb \psi_k\|_{L^2} \leq \frac{D}{4} \|\nabla\ptb \psi_k\|_{L^2}^2 +\frac{F_\infty^2}{D} \|\ptb \psi_k\|_{L^2}^2.
\end{equation}

Using again a combination of the H\"older and $\eps$-Young inequalities, we estimate the term on the right-hand side of the tested equation \eqref{Eq:Tested1} with the initials $\psi_{0k}$ by
$$\begin{aligned} &\ga(F\psi_{0k}-D\nabla \psi_{0k},\nabla \ptb \psi_k)_{L^2} \\ &\leq \frac{\ga}{2\eps} (F_\infty^2 \|\psi_{0k}\|_{L^2}^2+D^2\|\nabla \psi_{0k}\|_{L^2}^2)  + \eps_1\ga  \|\nabla \ptb \psi_k\|_{L^2}^2,\end{aligned}$$
where $\eps_1>0$ is a constant that we will determine accordingly below. We are not interested in tracking the constants $D$ and $F_\infty$ and therefore, we include them in a generic constant $C$ that may change from line to line. Moreover, we can estimate the norm of $\psi_{0k}$ by $\psi_0$ due to the projection property \eqref{Eq:Projection}. Consequently, we obtain the estimate
\begin{equation} \label{Eq:Tested1c}\begin{aligned} \ga(F\psi_0\!-\!D\nabla \psi_0,\nabla \ptb \psi_k)_{L^2} &\leq C\ga(t)  \|\psi_0\|_{H_0^1}^2 \!+\! \eps\ga \|\nabla \ptb \psi_k\|_{L^2}^2\end{aligned}\end{equation}
Lastly, we estimate the external force $f$ by
\begin{equation} \label{Eq:Tested1d}
\begin{aligned} \langle f,\ptb \psi_k\rangle_{H_0^1} &\leq \|f\|_{H^{-1}} \|\ptb \psi_k\|_{H_0^1} \\ &\leq C\|f\|_{H^{-1}}^2 + \frac{D}{4} \|\ptb \nabla \psi_k\|_{L^2}^2.
\end{aligned}
\end{equation}

Hence, we insert the estimates \eqref{Eq:Tested1a}--\eqref{Eq:Tested1d} in the tested equation \eqref{Eq:Tested1} to obtain the inequality
$$\begin{aligned} &\frac12 \pta \|\ptb \psi_k\|_{L^2}^2 + D \|\nabla \ptb \psi_k\|_{L^2}^2 \\  &\leq F_\infty \|\ptb \psi_k\|^2_{L^2}+C\|f\|_{H^{-1}}^2  + \frac{D}{2} \|\nabla \ptb \psi_k\|_{L^2}^2   \\&\quad + \eps_1 \ga(t) \|\nabla \ptb \psi_k\|_{L^2}^2  +Cg_{\alpha}(t)  \|\psi_0\|_{H_0^1}^2,
\end{aligned}$$
and we absorb the terms involving $D$ on the right-hand side by the respective term on the left-hand side, giving
$$\begin{aligned} &\frac12 \pta \|\ptb \psi_k\|_{L^2}^2 + \frac{D}{2} \|\nabla \ptb \psi_k\|_{L^2}^2 \\  &\leq F_\infty \|\ptb \psi_k\|^2_{L^2}+C\|f\|_{H^{-1}}^2    + \eps_1 \ga(t) \|\nabla \ptb \psi_k\|_{L^2}^2  +Cg_{\alpha}(t)  \|\psi_0\|_{H_0^1}^2.
\end{aligned}$$
We convolve this inequality with the kernel function $\ga$ to conclude
\begin{equation} \label{Eq:Galerkin2} \begin{aligned} &\frac12 \|\ptb \psi_k(t)\|_{L^2}^2 + \frac{D}{2} (\ga*\|\nabla \ptb \psi_k\|_{L^2}^2)(t) \\  &\leq F_\infty^2   (\ga *  \|\ptb \psi_k\|_{L^2}^2)(t)  + C(\ga*\|f\|_{H^{-1}}^2)(t)  \\ &\quad +\eps_1 (\ga*(\ga\cdot \|\nabla \ptb \psi_k\|_{L^2}^2))(t)  + C g_{2\alpha}(t)  \| \psi_0\|^2_{H^1_0}.
\end{aligned}
\end{equation}
where we used that $\ga*g_{\alpha}=g_{2\alpha}$, see \eqref{Eq:Semigroup}, and $$\begin{aligned}\big(\ga*(\pta \|\ptb \psi_k\|_{L^2}^2)\big)(t) &=\|\ptb \psi_k(t)\|_{L^2}^2 - \|(\ga*\pt \psi_k)(0)\|_{L^2}^2 \\ &= \|\ptb \psi_k(t)\|_{L^2}^2,
\end{aligned}$$ see \eqref{Eq:InverseConvolution}. Now, we observe that  the term $(\ga*\|f\|_{H^{-1}}^2)(t)$, $t \in (0,T_k)$,  can be bounded by
$$(\ga*\|f\|_{H^{-1}}^2)(t) \leq \sup_{t \in (0,T)} (\ga*\|f\|_{H^{-1}}^2)(t) =: \|f\|_{L^2_\alpha H^{-1}},$$
see again \eqref{Eq:LpAlpha} for the definition of the space $L^2_\alpha(0,T)$.

Further, we use \eqref{Eq:IneqGaG1Conv} to absorb the term involving $\eps_1$ on the right-hand side of the inequality \eqref{Eq:Galerkin2}. In fact, we absorb it by the term $\frac{D}{2} (\ga*\|\nabla \ptb \psi_k\|_{L^2}^2)(t)$ on the left-hand side  by noting that
$$\eps_1 (\ga * (\ga\|\nabla \psi_k\|_{L^2L^2}^2))(t)  \leq \eps_1 g_{\alpha+1}(T)  (\ga* \|\nabla \psi_k\|_{L^2}^2)(t).$$
We choose $\eps_1=\frac{D}{4g_{\alpha+1}(T)}$ to get 
$$\eps_1 (\ga * (\ga\|\nabla \psi_k\|_{L^2L^2}^2))(t)  \leq \frac{D}{4}  (\ga* \|\nabla \psi_k\|_{L^2}^2)(t),$$
and consequently, we obtain from \eqref{Eq:Galerkin2}  the inequality
$$\begin{aligned} &\frac12 \|\ptb \psi_k(t)\|_{L^2}^2 + \frac{D}{4} (\ga*\|\nabla \ptb \psi_k\|_{L^2}^2)(t) \\  &\leq F_\infty^2   (\ga *  \|\ptb \psi_k\|_{L^2}^2)(t)  + C\|f\|_{L^2_\alpha H^{-1}}^2   + C g_{2\alpha}(t)  \| \psi_0\|^2_{H^1_0}.
\end{aligned}$$

We notice that we are in the situation of the extended Henry--Gr\"onwall lemma, see Lemma \ref{Lem:GronFrac}, and we obtain the energy estimate
\begin{equation} \label{Eq:FinalEst1} \begin{aligned} &\frac12 \|\ptb \psi_k(t)\|_{L^2}^2 +  \frac{D}{4} \|\nabla \ptb \psi_k\|_{L^2_tL^2}^2 \\  &\leq C(F_\infty, \alpha, T) \cdot \Big((g_0+E)*\big( \|f\|_{L^2_\alpha H^{-1}}^2 +  g_{2\alpha}  \| \psi_0\|^2_{H^1_0} \big)\Big)(t) \\
		&=: \text{RHS}_{\eqref{Eq:FinalEst1}}(t)
\end{aligned}
\end{equation}
The estimate on the right-hand side is independent of $T_k$ and we infer from the no-blow-up theorem that we can continue the maximal time to $T$. However, since the right-hand side is "only" continuous in $t$ on $(0,T]$ and not at $t=0$ because of the presence of the term $g_{2\alpha}$, we are not able to take the essential supremum of the inequality  \eqref{Eq:FinalEst1} over $t \in (0,T)$. Therefore, we can obtain no bound of $\ptb \psi_k$ in $L^\infty$-in-time. Nonetheless,  $\ptb \psi_k$ is bounded in $L^2(0,T;H_0^1(\Omega))$ by inserting $t=T$ into the inequality \eqref{Eq:FinalEst1}. Moreover, we notice that $\|\ptb \psi_k(t)\|_{L^2(\Omega)}$ is bounded by the leading term $\sqrt{g_{2\alpha}(t)}=t^{\alpha-1/2}=g_{\alpha+1/2}$, which is continuous at $t=0$ for $\alpha > 1/2$ and in $L^p(0,T)$ for $\alpha+\frac12 > 1-\frac{1}{p}$, which is equivalent to $p<\frac{2}{1-2\alpha}$. Therefore, $\ptb \psi_k$ is bounded in the space $L^p(0,T;L^2(\Omega))$ with
\begin{equation} \label{Eq:LimitP} \begin{cases}
			p<\frac{2}{1-2\alpha}, &\alpha<\frac12, \\
	p<\infty, &\alpha=\frac12,\\
	p=\infty, &\alpha>\frac12.
\end{cases} \end{equation}

By the Eberlein--Smulian and Banach--Alaoglu theorems, see \cite[8.7]{alt2016linear}, these bounds yield the existence of a weakly converging subsequence $\ptb\psi_{k_j}$, i.e., it holds
\begin{equation}\label{Eq:Convergence1}\ptb \psi_{k_j} \longweak \zeta \quad \text{ in } L^p(0,T;L^2(\Omega)) \cap L^2(0,T;H_0^1(\Omega)),
\end{equation}
as $j \to \infty$.
We still need to figure out the representation of $\zeta$. If we are able to bound $\psi_k$ again in the Bochner space $L^2(0,T;H_0^1(\Omega))$, then we would obtain $\psi_k \rightharpoonup \psi$ for some limit function $\psi$, from which we can infer $\zeta=\ptb \psi$. We want to mention that we could also have obtained a bound of $\ptb \psi_k$ in the space $L^2_\alpha(0,T;H_0^1(\Omega))$. However, this space is not known to be reflexive and therefore, we cannot apply the Banach--Alaoglu theorem to infer a limit function in this space. \medskip

\noindent \textbf{(3) Energy estimates: Part 2.}
In order to obtain the desired bound of $\psi_k$, we test the Galerkin equation \eqref{Eq:FP_dis} with $\psi_k$, which yields
\begin{equation} \label{Eq:Galerkin3} \begin{aligned} &\frac12 \ddt \|\psi_k\|_{L^2}^2 + D (\ptb \nabla \psi_k,\nabla \psi_k)_{L^2}
\\ 
&= (F\cdot \nabla\ptb \psi_k,\psi_k)_{L^2}  +  \langle f,\psi_k \rangle_{H_0^1} - \ga( D\nabla \psi_{0k} - F\psi_{0k},\nabla \psi_k)_{L^2}.
\end{aligned}\end{equation}
For the term on the left-hand side involving the diffusion $D$, we apply Alikhanov's inequality \eqref{Eq:ChainOriginal} to infer
$$D (\ptb \nabla \psi_k,\nabla \psi_k)_{L^2} \geq \frac{D}{2} \ptb \|\nabla \psi_k\|_{L^2}^2.$$
Moreover, we apply again the H\"older inequality on the right-hand side of \eqref{Eq:Galerkin3} to obtain the energy estimate
$$\begin{aligned} &\frac12 \ddt \|\psi_k\|_{L^2}^2 + \frac{D}{2} \ptb \|\nabla \psi_k\|_{L^2}^2
\\ 
&\leq F_\infty  \|\nabla\ptb \psi_k\|_{L^2} \|\psi_k\|_{L^2}  +  \|f\|_{H^{-1}} \|\nabla \psi_k\|_{L^2}  \\ & \quad + \ga(t) \|\nabla \psi_k\|_{L^2} \big( D\|\nabla \psi_0\|_{L^2} + F_\infty \|\psi_0\|_{L^2} \big),
\end{aligned}$$
where we used the boundedness of the projection operator, see \eqref{Eq:Projection}.
Again, with the Young inequality, we obtain
$$\begin{aligned} &\frac12 \ddt \|\psi_k\|_{L^2}^2 + \frac{D}{2} \ptb \|\nabla \psi_k\|_{L^2}^2
\\ 
&\leq C F_\infty^2 \|\nabla \ptb \psi_k\|_{L^2}^2 +  \frac{\eps_2}{2}\|\nabla\psi_k\|_{L^2}  + C \|f\|_{H^{-1}}^2 + \frac{\eps_2}{2} \|\nabla \psi_k\|_{L^2}^2  \\ & \quad + \eps_3 \ga(t) \|\nabla \psi_k\|_{L^2}^2 + C\ga(t) \|\psi_0\|_{H_0^1}^2,
\end{aligned}$$
for some $\eps_2,\eps_3>0$ that we determine below.
After integrating this inequality over the time interval $(0,t)$, $t \leq T$, it yields
\begin{equation} \label{Eq:Estimate5} \begin{aligned} &\frac12 \|\psi_k(t)\|_{L^2}^2 + \frac{D}{2} (\ga* \|\nabla \psi_k\|_{L^2}^2)(t)
\\ 
&\leq \frac12 \|\psi_{0k}\|_{L^2}^2 + C F_\infty^2 \|\nabla \ptb \psi_k\|_{L^2_tL^2}^2   + C \|f\|_{L^2H^{-1}}^2 + \eps_2 \|\nabla \psi_k\|_{L^2_tL^2}^2   \\&\quad + \eps_3 \int_0^t \ga(s) \|\nabla \psi_k(s)\|_{L^2}^2 \ds + C g_{\alpha+1}(T) \|\psi_0\|_{H_0^1}^2 ,
\end{aligned}
\end{equation}
where we used that $\ga*1=g_{\alpha+1}$, see \eqref{Eq:Semigroup}, which is a continuous and bounded function on $[0,T]$ for any $\alpha>0$. 

We use the energy estimate \eqref{Eq:FinalEst1} from before to infer
\begin{equation} \label{Eq:Help1} 
\|\nabla \ptb \psi_k\|_{L^2_tL^2}^2 \leq \|\nabla \ptb \psi_k\|_{L^2L^2}^2 \leq \text{RHS}_{\eqref{Eq:FinalEst1}}(T).
\end{equation}
Furthermore, we use the auxiliary result \eqref{Eq:KernelNorm} to get
\begin{equation} \label{Eq:Help2} \begin{aligned} 
\eps_2 \|\nabla \psi_k\|_{L^2_tL^2}^2 &\leq \frac{\eps_2}{\ga(T)}  (\ga*\|\nabla \psi_k\|^2_{L^2})(t) \\ &\leq \frac{D}{8} (\ga*\|\nabla \psi_k\|^2_{L^2})(t),
\end{aligned}
\end{equation}
where we have chosen $\eps_2=\frac{D\ga(T)}{8}$.
Lastly, we use again \eqref{Eq:IneqGaG1} to infer
\begin{equation} \label{Eq:Help3} \begin{aligned} \eps_3 \int_0^t \ga(t) \|\nabla \psi_k\|_{L^2}^2 \ds \leq \eps_3 T (\ga* \|\nabla \psi_k\|_{L^2}^2)(t) \leq \frac{D}{8} (\ga* \|\nabla \psi_k\|_{L^2}^2)(t).
\end{aligned}
\end{equation}
Therefore, we set $\eps_3=\frac{D}{8T}$ and together with the auxiliary estimates \eqref{Eq:Help1}--\eqref{Eq:Help3} we obtain from \eqref{Eq:Estimate5}
\begin{equation} \label{Eq:FinalEst2} \begin{aligned} &\frac12 \|\psi_k(t)\|_{L^2}^2 + \frac{D}{4} (\ga* \|\nabla \psi_k\|_{L^2}^2)(t)
\\ 
&\leq  C \cdot \text{RHS}_{\eqref{Eq:FinalEst1}}(T)  + C \|f\|_{L^2H^{-1}}^2  + C g_{\alpha+1}(T) \|\psi_0\|_{H_0^1}^2 \\
&=:\text{RHS}_{\eqref{Eq:FinalEst2}}.
\end{aligned}\end{equation}

\noindent\textbf{(4) Weak and strong convergences.} From the estimate that we derived in \eqref{Eq:FinalEst2} we infer that that $\psi_k$ is bounded in the spaces $L^\infty(0,T;L^2(\Omega))$ and $L^2(0,T;H_0^1(\Omega))$, i.e, there is a limit function $\psi$ with
\begin{equation} \label{Eq:Convergence2a} \begin{aligned}
\psi_{k_j} &\longweak \psi \quad \text{ in } L^2(0,T;H_0^1(\Omega)), \\
\psi_{k_j} &\overset{*}{\longweak} \psi \quad \text{ in } L^\infty(0,T;L^2(\Omega)),
\end{aligned}\end{equation}
as $j \to \infty$.
By linearity of the differential operators, we obtain from \eqref{Eq:Convergence1}
\begin{equation} \label{Eq:Convergence2} \ptb \psi_{k_j} \longweak \ptb \psi \quad \text{ in } L^p(0,T;L^2(\Omega)) \cap L^2(0,T;H_0^1(\Omega)),\end{equation}
as $j \to \infty$, with $p$ as defined in \eqref{Eq:LimitP}. We note the compact embedding, see \eqref{Eq:aubinfractional},
 $$ L^2(0,T;H_0^1(\Omega)) \cap H^{1-\alpha}(0,T;H_0^1(\Omega)) \hookrightarrow\hookrightarrow L^2(0,T;L^2(\Omega)),$$
 from which we obtain the strong convergence 
   \begin{equation}\label{Eq:StrongConv1}\begin{aligned}
\psi_{k_j} &\longrightarrow \psi \quad &&\text{ in } L^2(0,T;L^2(\Omega).
  \end{aligned}\end{equation}

The derived convergences \eqref{Eq:Convergence2a}--\eqref{Eq:StrongConv1} are enough to pass to the limit in the Galerkin equation \eqref{Eq:FP_dis}. Nonetheless, we want to derive an additional estimate on $\pt \psi_k$ by testing the discretized Fokker--Planck equation \eqref{Eq:FP_dis} with $\Pi_{H_k} \zeta$ where $\zeta$ is an arbitrary element in $L^r(0,T;H_0^1(\Omega))$ with $r\geq 1$ depending on $\alpha$ (will be specified below). Then we obtain by the usual inequalities
$$\begin{aligned}(\pt\psi_k,\Pi_{H_k}\zeta )_{L^2L^2} 
  &=   (F \ptb \psi_k,\nabla \Pi_{H_k}\zeta)_{L^2L^2}-D(\ptb \nabla \psi_k,\nabla \Pi_{H_k}\zeta)_{L^2L^2}  \\ &\quad +\langle f,\Pi_{H_k}\zeta\rangle_{L^2H_0^1} - (  D \nabla \psi_{0k} - F \psi_{0k},\ga\nabla\Pi_{H_k}\zeta)_{L^2L^2}
  \\ 
  &\leq F_\infty \|\ptb \psi_k\|_{L^2L^2} \|\nabla \zeta\|_{L^2L^2} +D \|\ptb \nabla \psi_k\|_{L^2L^2}  \|\nabla \zeta\|_{L^2L^2}  \\ &\quad+\|f\|_{L^2H^{-1}}  \|\zeta\|_{L^2H_0^1}  + (F_\infty+D)\|\ga\|_{L^q}   \|\psi_0\|_{H^1}  \|\nabla\zeta\|_{L^{q'} L^2}   \\
  &\leq C \|\zeta\|_{L^r H_0^1},
  \end{aligned}$$ 
  where $r=\max\{q',2\}$ and $q'$ is the H\"older conjugate of $q=\frac{1}{1-\alpha}-\eps$ for $\eps\in (0,\frac{\alpha}{1-\alpha}]$.
  Therefore, $\pt \psi_k$ is bounded in $L^{r'}(0,T;H^{-1}(\Omega))$ where $r'$ is the H\"older conjugate of $r$. We note the compact embeddings, see \eqref{Eq:aubin}--\eqref{Eq:aubinfractional},
  $$\begin{aligned} 
  L^\infty(0,T;L^2(\Omega)) \cap W^{1,r'}(0,T;H^{-1}(\Omega)) &\hookrightarrow\hookrightarrow C([0,T];H^{-1}(\Omega)), \\
    W^{1-\alpha,2}(0,T;H_0^1(\Omega)) \cap W^{1,r'}(0,T;H^{-1}(\Omega))  &\hookrightarrow\hookrightarrow W^{1-\alpha,2}(0,T;L^2(\Omega)),
  \end{aligned}$$
  which provides us with the strong convergences (as $j\to \infty$)
  \begin{equation}\label{Eq:StrongConv}\begin{aligned}
\psi_{k_j} &\longrightarrow \psi \quad &&\text{ in } C([0,T];H^{-1}(\Omega), \\
\ptb \psi_{k_j} &\longrightarrow \ptb \psi \quad &&\text{ in } L^2(0,T;L^2(\Omega)). 
  \end{aligned}\end{equation}
\medskip

\noindent\textbf{(5) Limit process.}
In this step, we pass to the limit $j \to \infty$ in the  time-integrated $k_j$-th Galerkin system \eqref{Eq:FP_dis}. We use the derived convergences from the preceding result and show that the weak limit function $\psi$ satisfies the variational form of the time-fractional Fokker--Planck equation, i.e., $\psi$ is a weak solution in the sense of Definition \ref{Def:Weak}.

We consider the time-integrated $k_j$-th Galerkin system
$$\begin{aligned} 
& \int_0^T \Big( \langle \psi_{k_j}',\zeta \rangle_{H_0^1} +D(\ptb \nabla \psi_{k_j},\nabla \zeta)_{L^2} - (F \ptb \psi_{k_j},\nabla \zeta)_{L^2} \Big) \eta(t) \dt \\ &=\int_0^T \Big( \langle f,\zeta\rangle_{H_0^1} - \ga( D \nabla \psi_0 - F \psi_0,\nabla\zeta)_{L^2} \Big) \eta(t) \dt  \end{aligned}$$ for all 
$\zeta \in H_{k_j}$ and $\eta \in C_c^\infty(0,T)$.  Obviously, we are able to pass to the limit in all the terms thanks to the derived weak  convergences. E.g., we have
$$\begin{aligned}\int_0^T  (F\ptb \psi_{k_j},\nabla \zeta)_{L^2} \eta(t) \dt &\leq F_\infty \|\ptb \psi_{k_j} \|_{L^2L^2} \|\nabla\zeta\|_{L^2} \|\eta\|_{L^2}  \\ &\leq C \|\ptb\psi_{k_j} \|_{L^2H^1_0},
\end{aligned}$$
for all $\zeta \in H_{k_j}$, $\eta \in C_c^\infty(0,T)$. Since it holds the weak convergence 
$$\ptb\psi_{k_j} \rightharpoonup \ptb\psi_{k_j} \text{ in } L^2(0,T;H_0^1(\Omega)),$$ see \eqref{Eq:Convergence2}, it yields (as $j \to \infty$)
$$\int_0^T  (F \ptb \psi_{k_j},\nabla \zeta)_{L^2} \eta(t) \dt \longrightarrow \int_0^T (F \ptb \psi,\nabla \zeta)_{L^2} \eta(t) \dt,$$
for all $\zeta \in \cup_{j} H_{k_j}$. We observe that $\cup_j H_{k_j}$ is dense in $H_0^1(\Omega)$, which implies that the limit function $\psi$ indeed solves the variational form of the time-fractional Fokker--Planck equation. 

\medskip

\noindent\textbf{(6) Initial condition.}
  By the strong convergences, see \eqref{Eq:StrongConv}, we obtain at $t=0$ the convergence $\psi_{k_j}(0) \to \psi(0)$ in $H^{-1}(\Omega)$. However, it also holds $\psi_k(0) = \Pi_{H_k} \psi_0 \to \psi_0$ in $H_0^1(\Omega)$ as $j \to \infty$, from which we conclude $\psi(0) = \psi_0$ by the uniqueness of limits. Therefore, $\psi$ is a weak solution to the time-fractional Fokker--Planck equation in the sense of Definition \ref{Def:Weak}.  
  \medskip
  
 \noindent\textbf{(7) Uniqueness.} We consider two weak solutions $\psi_1$ and $\psi_2$ of the time-fractional Fokker--Planck equation in the sense of Definition \ref{Def:Weak}. Both solutions shall have the same initial data $\psi_0$ and outer force $f$. We subtract the variational forms of $\psi_1$ and $\psi_2$ from each other, and we define $\psi=\psi_1-\psi_2$, which satisfies 
 \begin{equation} \label{Eq:Uniqueness} \begin{aligned}
  &\langle \pt \psi,\zeta \rangle_{H_0^1}
  +  D(\ptb \nabla \psi,\nabla \zeta)_{L^2} - (F \ptb \psi,\nabla \zeta)_{L^2} = 0 \quad \forall \zeta \in H_0^1(\Omega).
\end{aligned} \end{equation}
We consider the test function $\zeta = \ptb\psi(t) \in H_0^1(\Omega)$ for a.e. $t \in (0,T)$, which yields together with Alikhanov's  inequality, see \eqref{Eq:ChainOriginal},
\begin{equation*} \begin{aligned}
  &\frac12  \pta \|\ptb \psi\|_{L^2}^2
  +  D \|\ptb\nabla \psi\|_{L^2}^2 \leq F_\infty \|\ptb \psi\|_{L^2} \|\nabla \ptb \psi\|_{L^2}.
\end{aligned}\end{equation*}
Furthermore, we apply Young's inequality to obtain 
\begin{equation} \label{Eq:Unique} \begin{aligned}
  &\frac12  \pta \|\ptb \psi\|_{L^2}^2
  +  \frac{D}{2} \|\ptb\nabla \psi\|_{L^2}^2 \leq \frac{F_\infty^2}{D} \|\ptb \psi\|_{L^2}^2.
\end{aligned}\end{equation}
After convolving this inequality with $\ga$ and applying the extended Henry--Gronwall inequality with $a\equiv 0$, see Lemma \ref{Lem:GronFrac}, the estimate \eqref{Eq:Unique} becomes
\begin{equation*} \begin{aligned}
  &\frac12 \|\ptb \psi(t)\|_{L^2}^2
  +  \frac{D}{2} \|\ptb\nabla \psi\|_{L^2_tL^2}^2 \leq 0.
\end{aligned}\end{equation*}
At this point, we further test the variational form \eqref{Eq:Uniqueness} by $\zeta=\psi(t) \in H_0^1(\Omega)$, which yields
 \begin{equation*} \begin{aligned}
  &\frac12 \ddt \|\psi\|_{L^2}^2
  +  \frac{D}{2} \ptb \|\nabla \psi\|_{L^2}^2 \leq F_\infty \|\ptb \psi\|_{L^2} \|\nabla \psi\|_{L^2} = 0.
\end{aligned}\end{equation*}
We integrate this inequality and observe that it holds $\|\psi(t)\|_{L^2}=0$ for any $t \in (0,T)$ i.e. $\psi_1=\psi_2$.
  \qed

\section{Numerical simulations} \label{Sec:Numerics}


Various numerical methods for time-fractional PDEs are summarized in the review article \cite{diethelm2020good} and in the monographs \cite{baleanu2012fractional,owolabi2019numerical,jin2023numerical}. 

\begin{figure}[htb!] \centering
	\begin{tikzpicture}[darkstyle/.style={circle,draw,fill=blue},redstyle/.style={circle,draw,fill=red}]
        \def\NUM{12}
        \def\FINAL{11}
        \def\GAMMA{1}
        \pgfmathtruncatemacro\NUMGAMMA{\NUM^\GAMMA}
        \pgfmathtruncatemacro\HALFGAMMA{2^\GAMMA}
		\draw[black, ultra thick] (0,0) -- (\FINAL,0);
		\node[redstyle,label=below: {$t_0\!=\!0$}]  at (0 ,0) {};
		\foreach \x in {1,...,\NUM}
		{\pgfmathtruncatemacro\xx{\x^\GAMMA}
			\node[darkstyle]   at ( \xx/\NUMGAMMA*\FINAL ,0) {};} 
		\node[redstyle,label=below: {$t_N\!=\!T$}]  at (\FINAL ,0) {};
		\node[redstyle,label=below: {$t_{N/2}$}]  at (\FINAL/\HALFGAMMA,0)  {};
	\end{tikzpicture} \\
 		\begin{tikzpicture}[darkstyle/.style={circle,draw,fill=blue},redstyle/.style={circle,draw,fill=red}]
        \def\NUM{12}
        \def\FINAL{11}
        \def\GAMMA{3}
        \pgfmathsetmacro\NUMGAMMA{4^\GAMMA}
        \pgfmathsetmacro\HALFGAMMA{sqrt(7)^\GAMMA}
		\draw[black, ultra thick] (0,0) -- (\FINAL,0);
		\node[redstyle,label=below: {$t_0\!=\!0$}]  at (0 ,0) {};
		\foreach \x in {1,...,\NUM}
		{\pgfmathsetmacro\xx{sqrt(\x)^\GAMMA}
			\node[darkstyle]   at ( \xx/\NUMGAMMA*\FINAL ,0) {};} 
		\node[redstyle,label=below: {$t_N\!=\!T$}]  at (\FINAL ,0) {};
		\node[redstyle,label=below: {$t_{N/2}$}]  at (\HALFGAMMA/\NUMGAMMA*\FINAL,0)  {};
	\end{tikzpicture} \\
 		\begin{tikzpicture}[darkstyle/.style={circle,draw,fill=blue},redstyle/.style={circle,draw,fill=red}]
        \def\NUM{12}
        \def\FINAL{11}
        \def\GAMMA{2}
        \pgfmathtruncatemacro\NUMGAMMA{\NUM^\GAMMA}
        \pgfmathtruncatemacro\HALFGAMMA{2^\GAMMA}
		\draw[black, ultra thick] (0,0) -- (\FINAL,0);
		\node[redstyle,label=below: {$t_0\!=\!0$}]  at (0 ,0) {};
		\foreach \x in {1,...,\NUM}
		{\pgfmathtruncatemacro\xx{\x^\GAMMA}
			\node[darkstyle]   at ( \xx/\NUMGAMMA*\FINAL ,0) {};} 
		\node[redstyle,label=below: {$t_N\!=\!T$}]  at (\FINAL ,0) {};
		\node[redstyle,label=below: {$t_{N/2}$}]  at (\FINAL/\HALFGAMMA,0)  {};
	\end{tikzpicture} \\
 	\begin{tikzpicture}[darkstyle/.style={circle,draw,fill=blue},redstyle/.style={circle,draw,fill=red}]
        \def\NUM{12}
        \def\FINAL{11}
        \def\GAMMA{3}
        \pgfmathtruncatemacro\NUMGAMMA{\NUM^\GAMMA}
        \pgfmathtruncatemacro\HALFGAMMA{2^\GAMMA}
		\draw[black, ultra thick] (0,0) -- (\FINAL,0);
		\node[redstyle,label=below: {$t_0\!=\!0$}]  at (0 ,0) {};
		\foreach \x in {1,...,\NUM}
		{\pgfmathtruncatemacro\xx{\x^\GAMMA}
			\node[darkstyle]   at ( \xx/\NUMGAMMA*\FINAL ,0) {};} 
		\node[redstyle,label=below: {$t_N\!=\!T$}]  at (\FINAL ,0) {};
		\node[redstyle,label=below: {$t_{N/2}$}]  at (\FINAL/\HALFGAMMA,0)  {};
	\end{tikzpicture} \\
\caption{Nonuniform time meshes on the interval $[0,T]$ with $t_n=(n/N)^\gamma T$ for $\gamma\in \{1,1.5,2,3\}$ (top to bottom) and $N=20$; the red nodes are $\{0,N/2,N\}$ in all cases.\label{Fig:Time}}
 \end{figure}
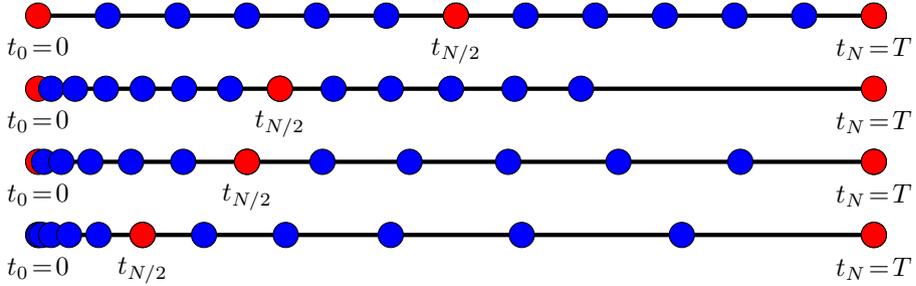

We assume a discretization $0=t_0<t_1<\dots<t_N=T$ of the time interval $[0,T]$. We do not utilize an equispaced time mesh, but a nonuniform one by discretizing the early times in finer steps. In particular, we assume that the $n$-th time step is of the form $t_n=(n/N)^\gamma T$ for $\gamma\geq 1$. If it holds $\gamma=1$, then we are again in a setting of a uniform mesh, see also Fig. \ref{Fig:Time} for a depiction of some time meshes for various values of $\gamma$.

We discretize the Caputo derivative by the nonuniform L1 scheme \cite[Section 3.2]{diethelm2020good}, i.e., it reads $$\ptb \psi \approx \frac{1}{\Gamma(1+\alpha)} \sum_{j=0}^{n-1} \omega_{n-j-1,n} (\psi_{n-j}-\psi_{n-j-1}),
$$ where $\psi_{n-j} \approx \psi(t_{n-j})$. The quadrature weights $\{\omega_{k,n}\}_{k=0}^{n-1}$ are given by the formula
$$\omega_{k,n}=\frac{(t_n-t_k)^{\alpha}-(t_n-t_{k+1})^{\alpha}}{\Delta t_{n-k}},$$
where we introduced the notation $\Delta t_{n-k}=t_{n-k}-t_{n-k-1}$. We use the finite element space $P_1$ for the space discretization and consequently, the fully discrete system reads
\begin{equation} \label{Eq:FP_Discretized} \begin{aligned}& \Big(\frac{\psi^n-\psi^{n-1}}{\Delta t_n},\zeta\Big)_H +  \sum_{j=0}^{n-1} \frac{\omega_{n-j-1,n}}{\Gamma(1+\alpha)}  (D\nabla(\psi_{n-j}-\psi_{n-j-1}),\nabla \zeta)_H \\&\quad - \sum_{j=0}^{n-1}   \frac{\omega_{n-j-1,n}}{\Gamma(1+\alpha)} (\psi_{n-j}-\psi_{n-j-1},F(t_n)\cdot \nabla \zeta)_H   \\ &= (f(t_n),\zeta)_H -  g_\alpha(t_n) \cdot (D(t_n)\nabla \psi_0,\nabla \zeta)_H + g_\alpha(t_n) \cdot  (\psi_0, F(t_n) \cdot\nabla \zeta)_H
\end{aligned} 
\end{equation}
for any test function $\zeta$. In particular, taking $\zeta=1$ gives
$$\begin{aligned}& \int_\Omega \psi^n \dx =\int_\Omega \psi^{n-1} \dx    + \Delta t_n \int_\Omega f(t_n) \dx,  
\end{aligned} $$
i.e., the Fokker--Planck setting with $f\equiv 0$ yields discrete mass conservation. We implement the discrete system in open-source computing platform \linebreak FEniCS, see \cite{alnaes2015fenics}.

 We consider the space interval $\Omega=(-5,15)$ with $\Delta x=1/1024$  and the time interval $[0,T]$ with $T=5$ where the $n$-th time step is given by $t_n=5(n/100)^2$. Moreover, we select as the initial data the Gaussian
$$\psi(0,x)=\psi_0(x)=\frac{1}{\sigma \sqrt{2\pi}} \text{exp}\Big(-\frac12 \Big( \frac{x-\mu}{\sigma} \Big)^2 \Big)$$
for $\sigma=0.1$ and $\mu=2$. 

Regarding model parameters, we choose $D=1$ and $f\equiv 0$. We take the space-time dependent force $F(t,x)=\sin(t)+x$ in Sec. \ref{Sec:Ex2} similar to \cite{angstmann2015generalized,mustapha2022second,le2016numerical,pinto2017numerical}. However, we first consider the case of an absent force $F \equiv 0$ in Sec. \ref{Sec:Ex1}, i.e., we are in the setting of a classical subdiffusion equation. In Sec. \ref{Sec:Ex3}, we consider the physically defeasible time-fractional Fokker--Planck equation with the Caputo derivative on the left-hand side, see again Sec. \ref{Sec:Derivation}, and compare this model numerically to the physically meaningful model that we have analyzed in this work.

\subsection{Example 1: Subdiffusion equation} \label{Sec:Ex1}

As we consider $F\equiv 0$ in this example, we essentially study the time-fractional heat equation
$$\pta \psi(x,t)=\Delta \psi(x,t),$$
which is also referred to as subdiffusion equation.

\begin{figure}[htb!]
    \centering
\subfigure[$t=0.02$]{
\includegraphics[height=.25\textheight]{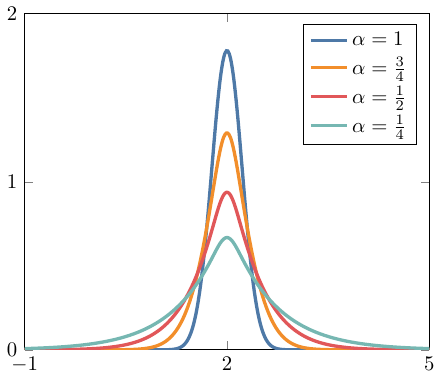} } \!\!\!\!\!\!\!\!
\subfigure[$t=0.5$]{
\includegraphics[height=.25\textheight]{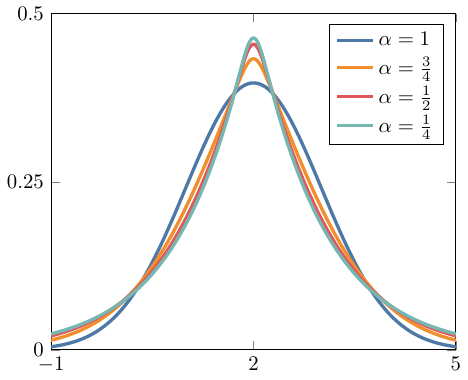} }
    \caption{Plot of the solution $\psi$ for varying values of $\alpha \in \big\{\frac14,\frac12,\frac34,1\big\}$; on the left (a) at $t=0.02$ and on the right (b) at $t=0.5$. \label{Fig:F0_AlphaVary}}
\end{figure}

We observe the typical behavior of a subdiffusive equation in the numerical simulations. At early times, the time-fractional model evolves faster stand the integer-order model. In Fig. \ref{Fig:F0_AlphaVary} (a), we see that the solution is more damped for $\alpha<1$ than for $\alpha=1$ at $t=0.02$. Moreover, the damping is larger for smaller values for $\alpha$. However, this behavior is exactly flipped if one considers a point further in time, e.g. $t=0.5$ as depicted in Fig. \ref{Fig:F0_AlphaVary} (b). After the initial fast evolution of the subdiffusion equation, the process is slower, and we observe that the smallest maximal value is represented by $\alpha=1$ at $t=0.5$. We can also observe that for $\alpha=1$ the typical round shape is present, whereas for $\alpha<1$ the tip at $x=2$ is less round.

\begin{figure}[htb!]
    \centering
\subfigure[$\alpha=1$]{\includegraphics[height=.25\textheight]{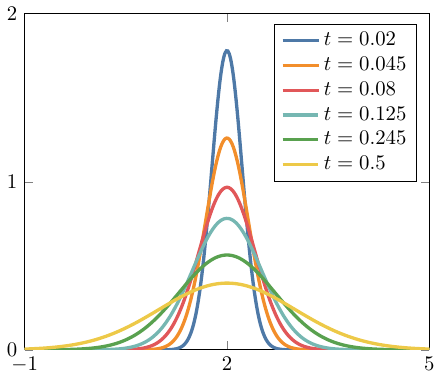} }
\!\!\!\!\!\!\!\!
\subfigure[$\alpha=\frac12$]{\includegraphics[height=.25\textheight]{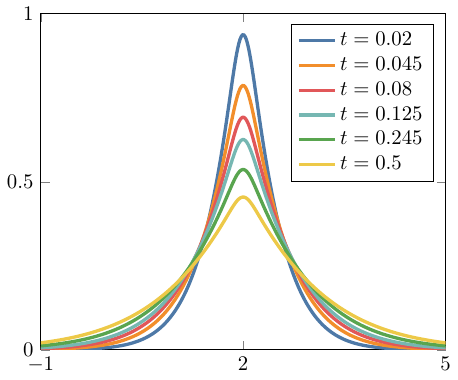}}
    \caption{Plot of the solution $\psi$ for varying time  $t \in \{0.02,0.045,0.08,0.125,0.245,0.5\}$; on the left (a) for $\alpha=1$ and on the right (b) for $\alpha=\frac12$. \label{Fig:F0_TimeVary}}
\end{figure}

\begin{figure}[htb!]
    \centering
\includegraphics{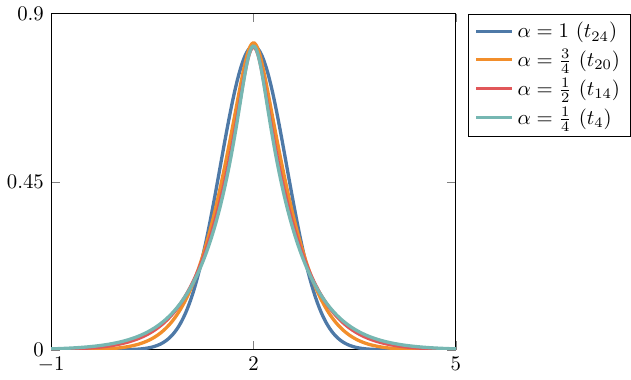}
    \caption{Plot of the solution $\psi$ for different values of $\alpha$ at different times $t$; we consider the pairings $(\alpha,t) \in \{(\frac14,t_4),(\frac12,t_{14}),(\frac34,t_{20}), (1,t_{24})\}$ for $t_n$ as defined at the beginning of Section \ref{Sec:Numerics}.}
    \label{Fig:F0_Fitting}
\end{figure}

We consider the time evolution for $\alpha=1$ in Fig. \ref{Fig:F0_TimeVary} (a) and for $\alpha=\frac12$ in Fig. \ref{Fig:F0_TimeVary} (b). The typical diffusion process can be observed and again, we notice the spikier tip for $\alpha=\frac12$. Moreover, the support of the function is larger for smaller $\alpha$.

Lastly, we try to fit the solution $\psi$ for different values of $\alpha$. The goal is to analyze whether it is necessary to consider the more complicated (analytically and numerically) time-fractional model, or whether this model's behavior can be replicated by an integer-order model.  This is done in Fig. \ref{Fig:F0_Fitting}, and we observe that the subdiffusive behavior cannot be imitated directly by the standard Fokker--Planck equation. Again, we observe the different support for each curve and the difference in the tip at $x=2$.

\subsection{Example 2: Space-time dependent force} \label{Sec:Ex2}
This time, we consider the space-time dependent force $F(x,t)=\sin(x)+t$ and therefore, we study the time-fractional Fokker--Planck equation
$$\begin{aligned}
&\pt \psi(x,t)-\Delta \ptb\psi(x,t)+ \div(F(x,t)\ptb \psi(x,t)) = \ga D\Delta \psi_0 - \ga \div(F\psi_0).
\end{aligned}$$

\begin{figure}[htb!]
    \centering
\subfigure[$t=0.02$]{
\includegraphics[height=.244\textheight]{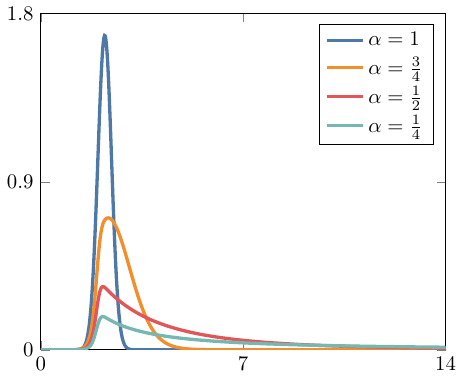} } \!\!\!\!\!\!\!\!\!
\subfigure[$t=0.18$]{
\includegraphics[height=.244\textheight]{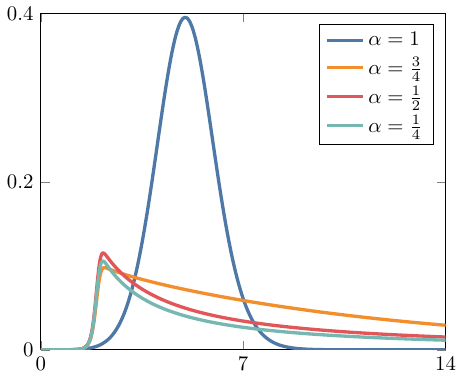} }
    \caption{Plot of the solution $\psi$ for varying values of $\alpha \in \big\{\frac14,\frac12,\frac34,1\big\}$; on the left (a) at $t=0.02$ and on the right (b) at $t=0.18$. \label{Fig:F1_AlphaVary}}
\end{figure}

Again, we observe the typical initial behavior of a subdiffusive equation. At the start, the time-fractional model evolves much faster stand the integer-order model. In Fig. \ref{Fig:F1_AlphaVary} (a), we see that the solution is more damped for $\alpha<1$ than for $\alpha=1$ at $t=0.02$. However, this time, we observe that the symmetry of the probability density functional $\psi$ is lost for $\alpha<1$. In the case of $\alpha=\frac12$ and $\alpha=\frac14$, the solution admits a large support up to the right end of the domain. that  In Fig. \ref{Fig:F1_AlphaVary} (b), we have plotted $\psi$ at a later time.  We observe that $\alpha=1$ is vastly different from the case of $\alpha<1$. This is also pronounced by the fact that $\ga(t) \to 0$ as $t \to \infty$ for $\alpha<1$, but in the case of $\alpha=1$ it holds $\ga(t) \equiv 1$, i.e., the right-hand side is just as large for all times.

\begin{figure}[htb!]
    \centering
\subfigure[$\alpha=1$]{\includegraphics[height=.25\textheight]{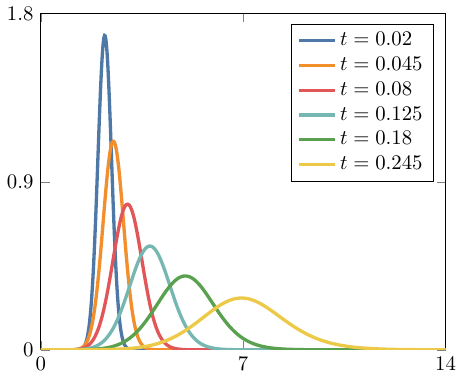} }
\!\!\!\!\!\!\!\!
\subfigure[$\alpha=\frac12$]{\includegraphics[height=.25\textheight]{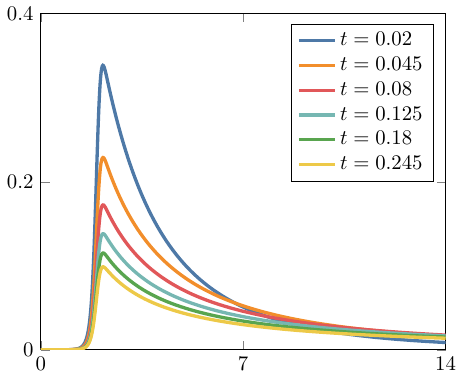}}
    \caption{Plot of the solution $\psi$ for varying time  $t \in \{0.02,0.045,0.08,0.125,0.245,0.5\}$; on the left (a) for $\alpha=1$ and on the right (b) for $\alpha=\frac12$. \label{Fig:F1_TimeVary}}
\end{figure}

We consider the time evolution for $\alpha=1$ in Fig. \ref{Fig:F0_TimeVary} (a) and for $\alpha=\frac12$ in Fig. \ref{Fig:F0_TimeVary} (b). The typical diffusion process can be observed and again, we notice the edgier tip for $\alpha=\frac12$. Moreover, the support of the function is larger for smaller $\alpha$.

Lastly, we try to fit the solution $\psi$ for different values of $\alpha$. This is done in Fig. \ref{Fig:F0_Fitting}, and we observe that the subdiffusive behavior cannot be imitated by an integer-order model. Again, we observe the different support for each curve and the difference in the tip at $x=2$.

\begin{figure}[htb!]
    \centering
\includegraphics{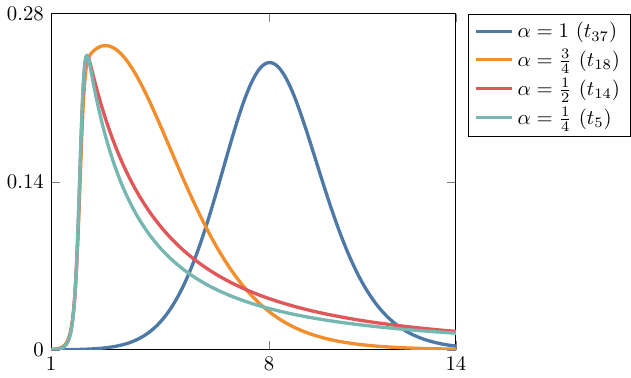}
    \caption{Plot of the solution $\psi$ for different values of $\alpha$ at different times $t$; we consider the pairings $(\alpha,t) \in \{(\frac14,t_5),(\frac12,t_{14}),(\frac34,t_{18}), (1,t_{37})\}$ for $t_n$ as defined at the beginning of Section \ref{Sec:Numerics}.}
    \label{Fig:F1_Fitting}
\end{figure}

\subsection{Example 3: Model comparison} \label{Sec:Ex3}

We consider the model as introduced in \eqref{Eq:ModelWrong} with no right-hand side, i.e.,
\begin{equation} \label{Eq:ModelWrong1} \begin{aligned}
&\pta \psi(x,t)-D \Delta \psi(x,t) +\div \big(F(t,x) \psi(x,t) \big) =0,
\end{aligned}\end{equation}
and we discretize it in the same manner as done for the time-fractional Fokker--Planck equation in \eqref{Eq:FP_Discretized}. Since this model has been studied in literature, we want to give it some attention by comparing it to the physically meaningful model. Again, we consider $F(x,t)=\sin(x)+t$.

\begin{figure}[H]
    \centering
\subfigure[$\alpha=\frac14$]{\includegraphics[height=.23\textheight]{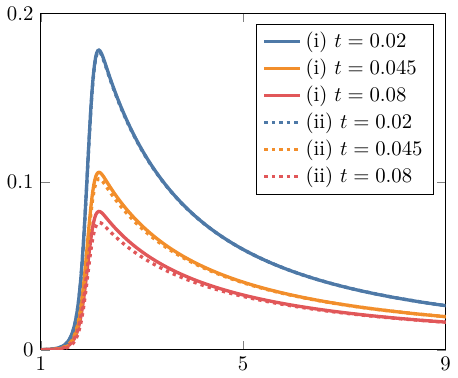} }
\!\!\!\!\!\!\!\!
\subfigure[$\alpha=\frac34$]{\includegraphics[height=.23\textheight]{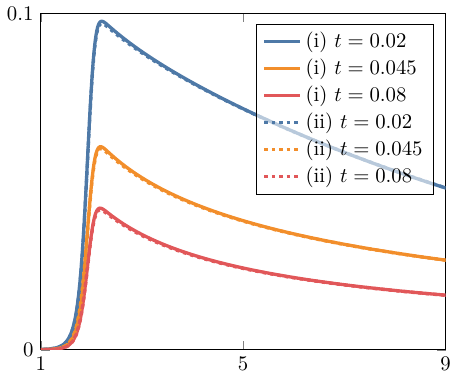}}
    \caption{Plot of the solution to the time-fractional Fokker--Planck equation \eqref{Eq:FP} (i) and the model with no right-hand side (ii) for varying time  $t \in \{0.02,0.045,0.08\}$; on the left (a) for $\alpha=\frac14$ and on the right (b) for $\alpha=\frac34$. \label{Fig:Wrong}}
\end{figure}

We compare it for $\alpha=\frac14$ in Fig. \ref{Fig:Wrong} (a) and for $\alpha=\frac34$ in Fig. \ref{Fig:Wrong} (b) for several time steps. We notice that the error gets larger for increasing time, and it is also more pronounced for smaller values for $\alpha$. We argue that this results from the fact that these models coincide for $\alpha=1$ and by continuity of the fractional parameter, the difference only gets larger the further one is from $\alpha=1$. Moreover, it holds $\ga(t) \to 0$ as $t \to \infty$ for $\alpha<1$ and therefore, it makes sense that asymptotically the right-hand side is negligible.

\section*{Acknowledgments}
Supported by the state of Upper Austria.



{\small	
	\bibliography{literature.bib}
	\bibliographystyle{AIMS} }


\end{document}